# CONSISTENCIES AND RATES OF CONVERGENCE OF JUMP-PENALIZED LEAST SQUARES ESTIMATORS

By Leif Boysen,[1] Angela Kempe,[2] Volkmar Liebscher,[3]
Axel Munk[4] and Olaf Wittich[3]

*Universität Göttingen, GSF–National Research Centre for Environment,
Universität Greifswald, Universität Göttingen and
Technische Universiteit Eindhoven*

We study the asymptotics for jump-penalized least squares regression aiming at approximating a regression function by piecewise constant functions. Besides conventional consistency and convergence rates of the estimates in $L^2([0,1))$ our results cover other metrics like Skorokhod metric on the space of càdlàg functions and uniform metrics on $C([0,1])$. We will show that these estimators are in an adaptive sense rate optimal over certain classes of "approximation spaces." Special cases are the class of functions of bounded variation (piecewise) Hölder continuous functions of order $0 < \alpha \leq 1$ and the class of step functions with a finite but arbitrary number of jumps. In the latter setting, we will also deduce the rates known from change-point analysis for detecting the jumps. Finally, the issue of fully automatic selection of the smoothing parameter is addressed.

**1. Introduction.** We consider regression models of the form

(1) $$Y_i^n = \overline{f}_i^n + \xi_i^n, \qquad i = 1, \ldots, n,$$

Received September 2006; revised September 2007.
[1]Supported by the Georg Lichtenberg program "Applied Statistics and Empirical Methods" and DFG Graduate Program 1023, "Identification in Mathematical Models."
[2]Supported by DFG, Priority Program 1114, "Mathematical Methods for Time Series Analysis and Digital Image Processing."
[3]Supported in part by DFG, Sonderforschungsbereich 386 "Statistical Analysis of Discrete Structures."
[4]Supported by DFG Grant "Statistical Inverse Problems under Qualitative Shape Constraints."

*AMS 2000 subject classifications.* Primary 62G05, 62G20; secondary 41A10, 41A25.
*Key words and phrases.* Jump detection, adaptive estimation, penalized maximum likelihood, approximation spaces, change-point analysis, multiscale resolution analysis, Potts functional, nonparametric regression, regressogram, Skorokhod topology, variable selection.

where $(\xi_i^n)_{n\in\mathbb{N}, 1\leq i\leq n}$ is a triangular scheme of independent zero-mean sub-Gaussian random variables and $\overline{f}_i^n$ is the mean value of a square integrable function $f \in L^2([0,1))$ over an appropriate interval $[x_{i-1}^n, x_i^n]$ [see, e.g., Donoho (1997)]

$$(2) \qquad \overline{f}_i^n = (x_i^n - x_{i-1}^n)^{-1} \int_{x_{i-1}^n}^{x_i^n} f(u)\,du.$$

For ease of notation, we will mostly suppress the dependency on $n$ in the sequel.

When trying to recover the characteristics of the regression function in applications, we frequently face situations where the most striking features are sharp transitions, called change points, edges or jumps [for data examples see Fredkin and Rice (1992), Christensen and Rudemo (1996), Braun, Braun and Müller (2000)]. To capture these features, in this paper we study a reconstruction of the original signal by step functions, which results from a least squares approximation of $Y = (Y_1, \ldots, Y_n)$ penalized by the number of jumps. More precisely, we consider minimizers $T_\gamma(Y) \in \arg\min H_\gamma(\cdot, Y)$ of the Potts functional

$$(3) \qquad H_\gamma(u, Y) = \frac{1}{n}\sum_{i=1}^n (u_i - Y_i)^2 + \gamma \cdot \#J(u).$$

Here $J(u) = \{i : 1 \leq i \leq n-1, u_i \neq u_{i+1}\}$ is the set of *jumps* of $u \in \mathbb{R}^n$. Note that the minimizer is not necessarily unique.

The name Potts functional refers to a model which is well known in statistical mechanics and was introduced by Potts (1952) as a generalization of the Ising model [Ising (1925)] for a binary spin system to more than two states. The original model was considered in the context of Gibbs fields with energy equal to the above penalty.

Various other strategies dealing with discontinuities are known in the literature. Kernel regression as (linear) nonparametric method offers various ways to identify jumps in the regression function, essentially by estimating modes of the derivative; see, for example, Hall and Titterington (1992), Loader (1996), Müller (1992) or Müller and Stadtmüller (1999). Other approaches like local $M$-smoothers [Chu et al. (1998)], sigma-filter [Godtliebsen, Spjøtvoll and Marron (1997)], chains of sigma-filters [Aurich and Weule (1995)] or adaptive weights smoothing [Spokoiny (1998), Polzehl and Spokoiny (2003)] are based on nonlinear averages which mimic robust $W$-estimators [cf. Hampel et al. (1986)] near discontinuities. Therefore, they do not blur the jump as much as linear methods would do.

The case when the regression function is a step function has been studied first by Hinkley (1970) and later by Yao (1988) and Yao and Au (1989). Given a known upper bound for the number of jumps, Yao and Au (1989) derive the optimal $O(n^{-1/2})$ and $O(n^{-1})$ rates for recovering the function in an $L^2$ sense and detecting the jump points, respectively. Their results have been generalized to overdispersion models and applied to DNA-segmentation by Braun, Braun and Müller (2000). Without the constraint of a known upper bound for the number of jumps, Birgé and Massart (2007) give a nonasymptotic bound for the MSE for a slightly different penalty.

In this more general setting we will deduce the same (parametric) rates as Yao and Au (1989) for the Potts minimizer if $f$ is piecewise constant with a finite but arbitrary number of jumps. We show that the estimate asymptotically reconstructs the correct number of jumps with probability 1. Further we will give (optimal) rates in the Skorokhod topology, which provides *simultaneous* convergence of the jump points and the graph of the function, respectively. As far as we know, this approach is new to regression analysis.

If the true regression function is not a step function, the Potts minimizer cannot compete in terms of rate of convergence for smoothness assumptions stronger than $C^1$. This is due to the nonsmooth approach of approximation via step functions and could be improved by fitting polynomials between estimated jumps [see Spokoiny (1998), Kohler (1999)]. For less smooth functions, however, we will show that it is adaptive and obtains optimal rates of convergence. To this end, we prove rates of convergence in certain classes of "approximation spaces" well known in approximation theory [DeVore and Lorentz (1993)]. To our knowledge, these spaces have not been introduced to statistics before. As special cases, we obtain (up to a logarithmic factor) the optimal $O(n^{-1/3})$ and $O(n^{-\alpha/(2\alpha+1)})$ rates if $f$ is of bounded variation or if $f$ is (piecewise) Hölder continuous on $[0,1]$ of order $1 \geq \alpha > 0$, respectively. The logarithmic factor occurs, since we give almost sure bounds instead of the more commonly used stochastic or mean square error bounds. Optimality in the class of functions with bounded variation shows that the Potts minimizer has the attribute of "local adaptivity" [Donoho et al. (1995)]. Under the assumption that the error is bounded, Kohler (1999) obtained nearly the same rates (worse by an additional logarithmic term) in these Hölder classes for the mean square error of a similar estimator.

We stress that minimizing $H_\gamma$ in (3) results in a step function, that is, a regressogram in the sense of Tukey (1961). Hence, this paper also answers the question how to choose the partition of the regressogram in an asymptotic optimal way [cf. Eubank (1999)] over a large scale of approximation spaces.

*Subset selection and TV penalization.* Our results can be viewed as a result on subset selection in a linear model $Y = \alpha + \beta^T X + \varepsilon$ with covariates

$X$. In this context our estimator minimizes the functional

$$L_n(\alpha, \beta) := \sum_{i=1}^{n}\left(Y_i - \alpha - \sum_{j=1}^{k}\beta_j X_{ij}\right)^2 \quad \text{subject to} \quad \#\{j : \beta_j \neq 0\} \leq N,$$

or (for proper $N$), what is equivalent for a proper choice of $\gamma$, minimization of

$$L_n(\alpha, \beta) + \gamma \#\{j : \beta_j \neq 0\}.$$

Setting $k = n - 1$ as well as $X_{ij} = 1$ for $j < i$ and 0 else, we obtain the Potts functional (3) with $u_1 = \alpha$ and $u_i = \alpha + \sum_{j=1}^{i-1} \beta_j$ for $2 \leq i \leq n$. In general, to select the correct variables, one requires a kind of oversmoothing, which is reflected by our results in the present paper. The Potts smoother in (3) achieves this by means of an $\ell_0$ penalty and for nearly uncorrelated predictors it is well known that $\ell_1$ penalization has almost the same properties as complexity-penalized least squares regression [cf. Donoho (2006a, 2006b)]. However, as a variable selection problem, detection of jumps in regression has a special feature, namely, the covariates $X_{ij}$ are highly correlated and these results do not apply. A similar comment applies to TV penalized estimation, as, for example, considered by Mammen and van de Geer (1997) which aims for minimizing

$$F_\gamma(u, Y) = \gamma \cdot \sum_{1 \leq i \leq n-1} |u_i - u_{i+1}| + \sum_{i=1}^{n}(u_i - Y_i)^2.$$

This can also be viewed in this context. Choosing $X_{ik}$ as above, it is a special case of the lasso, which was introduced by Tibshirani (1996) and minimizes $L_n(\alpha, \beta)$ subject to $\sum_{j=1}^{k}|\beta_j| \leq t$. Again, for (nearly) uncorrelated predictors, the lasso comes close to the $\ell_0$ solution. Thus, the relation of the Potts functional to the total variation penalty is roughly the same as the relation of subset selection to the lasso. In fact, for highly correlated predictors, the relationship between $\ell_0$ and $\ell_1$ solutions is much less understood and this question is above the scope of the paper. However, it seems that in our case $\ell_1$ penalization performs suboptimally. As an indication, from Mammen and van de Geer (1997), Theorem 10, we obtain an upper rate bound of $O_\mathbb{P}(n^{-\alpha/3})$ for the error of the total variation penalized least squares estimator of an $\alpha$-Hölder continuous function in contrast to the (optimal) rate of $O_\mathbb{P}(n^{\alpha/(2\alpha+1)})$, achieved by the Potts minimizer.

A reason for this difference is that the Potts functional will generally lead to fewer but higher jumps in the reconstruction, and hence is even more sparse than $\ell_1$ or TV based reconstructions. In general, a side phenomenon related to such sparsity of an estimator is a bad uniform risk behavior [see Pötscher and Leeb (2008)]. Although the conditions of that paper are not fulfilled in our model (basically, contiguity of the error distributions will fail), this phenomenon can be observed numerically in our situation. Our estimate will fail when the number of jumps grows too fast with the number of observations and small plateaus in the data will not be captured. However, our emphasis is on estimation of the *main* data features (here jumps) to obtain a sparse description of data, similar in spirit to Davies and Kovac (2001).

*Computational issues.* In general, a major burden of $\ell_0$ penalization is that it leads to optimization problems which are often NP hard and relaxation of this functional becomes necessary or other penalties, such as $\ell_1$, have to be used. Interestingly, computation of the minimizer of the Potts functional in (3) is a notable exception. The family $(T_\gamma(Y))_{\gamma>0}$ can be computed in $O(n^3)$ and the minimizer for one $\gamma$ in $O(n^2)$ steps [see Winkler and Liebscher (2002)]. At the heart of that result is the observation that the set of partitions of a discrete interval carries the structure of a directed acyclic graph which makes dynamic programming directly applicable [see Friedrich et al. (2008)].

The paper is organized as follows: after introducing some notation in Section 2, we provide in Section 3.1 the rates and consistency results for step functions and general bounded functions in the $L^2$ metric. In Section 3.2 we present the results of convergence in Hausdorff metric for the set of jump functions and in Section 3.3 for the Skorokhod topology for the regression function. In Section 3.4 we will introduce a simple data-driven parameter selection strategy resulting from our previous results and compare this to a multiresolution approach as in Davies and Kovac (2001). We briefly discuss relations to other models such as Bayesian imaging and extensions to higher dimensions in Section 4. Technical proofs are given in the Appendix.

This paper is complemented by the work of Boysen et al. (2007) which contains technical details of some of the proofs, the consistency of the estimates for more general noise conditions and the consistency of the empirical scale space $(T_\gamma(Y))_{\gamma>0}$ toward its deterministic target [cf. Chaudhuri and Marron (2000)].

**2. Model and notation.** For a functional $F:\Theta \to \mathbb{R} \cup \{\infty\}$, we denote by $\arg\min F$ the subset of $\Theta$ consisting of all minimizers of $F$. Let $S([0,1)) = \{f : f = \sum_{i=1}^n \alpha_i \mathbf{1}_{[t_i,t_{i+1})}, \alpha_i \in \mathbb{R}, 0 = t_1 < \cdots < t_{n+1} = 1, n \in \mathbb{N}\}$ denote the space of right-continuous step functions and let $D([0,1))$ denote the càdlàg space of right-continuous functions on $[0,1]$ with left limits and left-continuous at 1. Both will be considered as subspaces of $L^2([0,1))$ with the obvious identification of a function with its equivalence class, which is injective for these two spaces. $\|\cdot\|$ will denote the norm of $L^2([0,1))$ and the norm on $L^\infty([0,1))$ is denoted by $\|\cdot\|_\infty$.

Minimizers of the Potts functionals (3) will be embedded into $L^2([0,1))$ by the map $\iota^n : \mathbb{R}^n \longmapsto L^2([0,1))$,

$$(4) \qquad \iota^n((u_1,\ldots,u_n)) = \sum_{i=1}^{n} u_i \mathbf{1}_{[(i-1)/n, i/n)}.$$

Under the regression model (1), this leads to estimates $\hat{f}_n = \iota^n(T_{\gamma_n}(Y))$, that is,

$$(5) \qquad \hat{f}_n \in \iota^n(\arg\min H_{\gamma_n}(\cdot, Y)).$$

Here and in the following $(\gamma_n)_{n \in \mathbb{N}}$ is a (possibly random) sequence of smoothing parameters. We suppress the dependence of $\hat{f}_n$ on $\gamma_n$ since this choice will be clear from the context.

For the noise, we assume the following uniform sub-Gaussian condition. For a discussion on how this condition can be weakened [see Boysen et al. (2007)].

CONDITION (A). The triangular array $(\xi_i^n)_{n \in \mathbb{N}, 1 \leq i \leq n}$ of random variables obeys the following properties.

(i) For all $n \in \mathbb{N}$ the random variables $(\xi_i^n)_{1 \leq i \leq n}$ are independent.

(ii) There is a universal constant $\beta \in \mathbb{R}$ such that $\mathbb{E}e^{\nu \xi_i^n} \leq e^{\beta \nu^2}$ for all $\nu \in \mathbb{R}$, $1 \leq i \leq n$, and $n \in \mathbb{N}$.

Finally, we recall the definition of Hölder classes. We say that a function $f : [0,1] \to \mathbb{R}$ belongs to the Hölder class of order $0 < \alpha \leq 1$, if there exists $C > 0$ such that

$$|f(x) - f(y)| \leq C|x-y|^\alpha \qquad \text{for all } x, y \in [0,1].$$

**3. Consistency and rates.** In order to extend the Potts functional in (3) to $L^2([0,1))$, we define for $\gamma > 0$, the continuous Potts functionals $H_\gamma^\infty : L^2([0,1)) \times L^2([0,1)) \to \mathbb{R} \cup \{\infty\}$:

$$H_\gamma^\infty(g, f) = \begin{cases} \gamma \cdot \#J(g) + \|f - g\|^2, & \text{if } g \in S([0,1)), \\ \infty, & \text{otherwise.} \end{cases}$$

Here $J(g) = \{t \in (0,1) : g(t-) \neq g(t+)\}$ is the set of jumps of $g \in S([0,1))$. By definition, we have for every $g \in \arg\min H_\gamma^\infty(\cdot, f)$ that $H_\gamma^\infty(g, f) \leq H_\gamma^\infty(0, f) = \|f\|^2$ and therefore $\#J(g) \leq \gamma^{-1}\|f\|^2$ for $\gamma > 0$. Since a minimizer is uniquely determined by its set of jumps, minimizing $H_\gamma^\infty$ can be reduced to a minimization problem on the compact set of jump configurations with not more than $\gamma^{-1}\|f\|^2$ jumps which implies existence of a minimizer. For $\gamma = 0$, we set $H_0^\infty(g, f) = \|f - g\|^2$ for all $g \in L^2([0,1))$, hence

LEMMA 1. *For any $f \in L^2([0,1))$ and all $\gamma \geq 0$ we have $\arg\min H_\gamma^\infty(\cdot, f) \neq \varnothing$.*

In order to keep the presentation simple, we choose throughout the following an equidistant design $x_i^n = i/n$ in the model (1) and (2). All results given remain valid for designs with design density $h$, such that $\inf_{t \in [0,1]} h(t) > 0$ and $h$ is Hölder continuous on $[0,1]$ of order $\alpha > 1/2$. Moreover, for all theorems in this section we will assume that $Y^n$ is determined through (1) and the noise $\xi^n$ satisfies Condition (A).

3.1. *Convergence in $L^2$.* We investigate the asymptotic behavior of the Potts minimizer when the sequence $(\gamma_n)_{n \in \mathbb{N}}$ converges to a constant $\gamma$ for $\gamma > 0$ and $\gamma = 0$, respectively. If $\gamma > 0$, we do not recover the original function in the limit, but a parsimonious representation at a certain scale of interest determined by $\gamma$. For $\gamma = 0$ the Potts minimizer is consistent for the true signal under some conditions on the sequence $(\gamma_n)_{n \in \mathbb{N}}$:

(H1) $(\gamma_n)_{n \in \mathbb{N}}$ satisfies $\gamma_n \to 0$ and $\gamma_n n / \log n \to \infty$ $\mathbb{P}$-a.s.

For the consistency in approximation spaces in Theorem 2, we consider instead

(H2) $(\gamma_n)_{n \in \mathbb{N}}$ satisfies $\gamma_n \to 0$ and $\gamma_n \geq (1+\delta) 12\beta \, \log n / n$ $\mathbb{P}$-a.s. for almost every $n$ and some $\delta > 0$. Here $\beta$ is given by the noise Condition (A).

THEOREM 1. (i) *Assume that $f \in L^2([0,1))$ and $\gamma > 0$ are such that $f_\gamma$ is a unique minimizer of $H_\gamma^\infty(\cdot, f)$. Moreover, suppose $(\gamma_n)_{n \in \mathbb{N}}$ satisfies $\gamma_n \to \gamma$ $\mathbb{P}$-a.s.; then*

$$\hat{f}_n \xrightarrow[n \to \infty]{L^2([0,1))} f_\gamma \qquad \mathbb{P}\text{-}a.s.$$

(ii) *Let $f \in L^2([0,1))$ and $(\gamma_n)_{n \in \mathbb{N}}$ fulfill* (H1). *Then*

$$\hat{f}_n \xrightarrow[n \to \infty]{L^2([0,1))} f \qquad \mathbb{P}\text{-}a.s.$$

(iii) *Let $f \in S([0,1))$ and $(\gamma_n)_{n \in \mathbb{N}}$ fulfill* (H1). *Then*

$$\|\hat{f}_n - f\| = O\left(\sqrt{\frac{\log n}{n}}\right) \qquad \mathbb{P}\text{-}a.s.$$

*Moreover,*

$$\|\hat{f}_n - f\| = O_\mathbb{P}\left(\sqrt{\frac{1}{n}}\right).$$

We stress that the parametric rates in Theorem 1(iii) are obtained for a broad range of rates for the sequence of smoothing parameters. It is only required that $\gamma_n$ converges to zero slower than $\log n/n$. When trying to extend these results to more general function spaces, the question arises, which properties of the true regression function $f$ determine the almost sure rate of convergence of the Potts estimator. It turns out that the answer lies in the speed of approximation of $f$ by step functions. Let us introduce the approximation error

$$\Delta_k(f) := \inf\{\|g - f\| : g \in S([0,1)), \#J(g) \leq k\} \tag{6}$$

and the corresponding approximation spaces

$$\mathcal{A}^\alpha = \left\{f \in L^\infty[0,1] : \sup_{k \geq 1} k^\alpha \Delta_k(f) < \infty\right\}$$

for $\alpha > 0$. The following theorem gives the almost sure rates of convergence for these spaces.

THEOREM 2. *If $f \in \mathcal{A}^\alpha$ and $(\gamma_n)_{n \in \mathbb{N}}$ satisfies condition* (H2), *then*

$$\|\hat{f}_n - f\| = O(\gamma_n^{\alpha/(2\alpha+1)}) \qquad \mathbb{P}\text{-}a.s.$$

Now we give examples of well known function spaces contained in $\mathcal{A}^\alpha$ for $\alpha \leq 1$.

EXAMPLE 1. Suppose $f$ has finite total variation. Then, $f \in \mathcal{A}^1$ holds. Choosing $\gamma_n \asymp \log n/n$ such that condition (H2) is fulfilled yields $\|\hat{f}_n - f\| = O((\log n/n)^{1/3})$ $\mathbb{P}$-a.s.

PROOF. For the application of Theorem 2 we need to show that there is a $\delta > 0$ such that for all $k \in \mathbb{N}$, $k \geq 1$, there is an $f_k \in S([0,1))$ with $\|f - f_k\| \leq \delta/(k+1)$ and $\#J(f_k) \leq k$. Since each function of finite total variation is the difference of two increasing functions and $\#J(g + g') \leq \#J(g) + \#J(g')$, it is enough to consider increasing $f$ with $f(0) = 0$ and $f(1) < 1$. Define for $i = 1, \ldots, k$ intervals

$$I_i = f^{-1}([(i-1)/k, i/k)).$$

Then, $f_k(x) = \sum_{i=1}^k \mathbf{1}_{I_i}(x)(i-1/2)/k$ satisfies $\|f - f_k\| \leq \|f - f_k\|_\infty \leq (2k)^{-1}$ which completes the proof. □

EXAMPLE 2. Suppose $f$ belongs to a Hölder class of order $\alpha$ (with $0 < \alpha \leq 1$). Then, $f \in \mathcal{A}^\alpha$ holds. For $\gamma_n \asymp \log n/n$ fulfilling condition (H2), we get that $\|\hat{f}_n - f\| = O((\log n/n)^{\alpha/(2\alpha+1)})$ $\mathbb{P}$-a.s.

PROOF. Analogous to the proof above, we define for $I_i = [(i-1)/k, i/k)$ the function $f_k(x) = \sum_{i=1}^{k} \mathbf{1}_{I_i}(x) f((i-1/2)/k)$. On $I_i$ we have $\|f(x) - f(y)\|_\infty \leq Ck^{-\alpha}$. Thus $\|f - f_k\| \leq \|f - f_k\|_\infty \leq C(2k)^{-\alpha}$ holds. □

Obviously this result still holds, if the regression function $f$ is piecewise Hölder with finitely many jumps.

REMARK 1 (The case $\alpha > 1$). The characterization of the sets $\mathcal{A}^\alpha$ and related questions are a prominent theme in nonlinear approximation theory [see, e.g., DeVore (1998), DeVore and Lorentz (1993)]. For $f$ piecewise $C^1$, it is known that $\alpha > 1$ implies that $f$ is piecewise constant [Burchard and Hale (1975)], whereas this is still an open problem for general $f$. We conjecture that this implication holds for any $f$. This would imply that stronger smoothness assumptions than in the examples above do not yield better convergence rates.

Choosing $\gamma_n$ independently of the function and the function class as in the examples above yields convergence rates which are up to a logarithmic factor the optimal rates in the classes $\mathcal{A}^\alpha$, $0 < \alpha \leq 1$ and $S([0,1))$. This shows that the estimate is adaptive over these classes. The additional logarithmic factor originates from giving almost sure rates of convergence.

3.2. *Hausdorff convergence of the jump-sets.* In this section we present the rates known from change-point analysis for detecting the locations of jumps if $f$ is a step function. Moreover, the following theorem shows that we will eventually estimate the right number of jumps almost surely. Before stating the results, we recall the definition of the Hausdorff metric $\rho_H$ on the space of closed subsets contained in $(0,1)$. For nonempty closed sets $A, B \subset (0,1)$ set

$$\rho_H(A, B) = \max\left\{\max_{a \in A} \min_{b \in B} |b - a|, \max_{b \in B} \min_{a \in A} |b - a|\right\}$$

and $\rho_H(A, \varnothing) = \rho_H(\varnothing, A) = 1$.

THEOREM 3. *Let $f \in S([0,1))$ and $(\gamma_n)_{n \in \mathbb{N}}$ fulfill* (H1). *Then:*

(i) $\#J(\hat{f}_n) = \#J(f)$ *for large enough $n$ $\mathbb{P}$-a.s.,*
(ii) $\rho_H(J(\hat{f}_n), J(f)) = O(\log n/n)$ $\mathbb{P}$-*a.s.,*
(iii) $\rho_H(J(\hat{f}_n), J(f)) = O_\mathbb{P}(1/n)$.

REMARK 2 (Distribution of the jump locations and estimated function values). With the help of Theorem 3(i) we can derive the asymptotic distribution of the jump locations and of the estimated function values between, obtaining the same results as Yao and Au (1989), who assumed an a priori bound of the number of jumps. To this end, note that the estimator of Yao and Au (1989) and the Potts minimizer coincide if they have the same number of jumps. Denoting the ordered jumps of $f$ and their estimators by $(\tau_1, \ldots, \tau_R)$ and $(\hat{\tau}_1, \ldots, \hat{\tau}_{\hat{R}})$, respectively, we know by Theorem 3(i) that asymptotically $\hat{R} = R$ holds almost surely. For $\hat{R} = R$ we get that $n(\hat{\tau}_1, \ldots, \hat{\tau}_R)$ are asymptotically independent and the limit distribution of $n(\hat{\tau}_r - [\tau_r])$ is the minimum of a two-sided asymmetric random walk [cf. Yao and Au (1989), Theorem 1]. Moreover, the estimated function values are asymptotically normal with the parametric $\sqrt{n}$-rate.

3.3. *Convergence in Skorokhod topology.* Now that we have established rates of convergence for the graph of the function as well as for the set of jump points, it is natural to ask whether one can handle both simultaneously. To this end, we recall the definition of the Skorokhod metric [Billingsley (1968), Chapter 3]. Let $\Lambda_1$ denote the set of all strictly increasing continuous functions $\lambda : [0,1] \longmapsto [0,1]$ which are onto. We define for $f, g \in D([0,1))$

$$\rho_S(f,g) = \inf \left\{ \max\left( L(\lambda), \sup_{0 \leq t \leq 1} |f(\lambda(t)) - g(t)| \right) : \lambda \in \Lambda_1 \right\},$$

where $L(\lambda) = \sup_{s \neq t \geq 0} |\log \frac{\lambda(t) - \lambda(s)}{t-s}|$. The topology induced by this metric is called $J_1$-topology.

We find that in the situation of Theorem 1(i) we can establish consistency without further assumptions, whereas in the situation of Theorem 1(ii), $f$ has to belong to $D([0,1))$.

THEOREM 4. (i) *Under the assumptions of Theorem* 1(i) *we have*

$$\hat{f}_n \xrightarrow[n \to \infty]{D([0,1))} f_\gamma \qquad \mathbb{P}\text{-}a.s.$$

(ii) *If $f \in D([0,1))$ and $(\gamma_n)_{n \in \mathbb{N}}$ satisfies condition* (H1), *then*

$$\hat{f}_n \xrightarrow[n \to \infty]{D([0,1))} f \qquad \mathbb{P}\text{-}a.s.$$

*If $f$ is continuous on $[0,1]$*

$$\hat{f}_n \xrightarrow[n \to \infty]{L^\infty([0,1])} f \qquad \mathbb{P}\text{-}a.s.$$

(iii) *If $f \in S([0,1))$ and $(\gamma_n)_{n \in \mathbb{N}}$ satisfies condition* (H1), *then*

$$\rho_S(\hat{f}_n, f) = O\left(\sqrt{\frac{\log n}{n}}\right) \qquad \mathbb{P}\text{-}a.s.$$

*Moreover,*

$$\rho_S(\hat{f}_n, f) = O_\mathbb{P}\left(\sqrt{\frac{1}{n}}\right).$$

3.4. *Parameter choice and simulated data.* In this section we assume $\xi_i^n \sim N(0, \sigma^2), i = 1, \ldots, n$ i.i.d. for all $n$. Note that in this case we have $\beta = \sigma^2/2$ in Condition (A). Theorem 2 directly yields a simple data-driven procedure for choosing the parameter $\gamma$ which leads to optimal rates of convergence. For a strongly consistent estimate $\hat{\sigma}$ of $\sigma$, the choice $\hat{\gamma}_n = C\hat{\sigma}^2 \log n/n$ almost surely satisfies condition (H2) for $C > 6$ and gives the rates of Theorem 2. However, in simulations it turns out that smaller choices of $C$ lead to better reconstructions. A closer look at the proof of Theorem 2 shows that the constant in condition (H2) mainly depends on the behavior of the maximum of the partial sum process $\sup_{1 \leq i \leq j \leq n} (\xi_i^n + \cdots + \xi_j^n)^2/(j-i+1)$. As we consider a triangular scheme instead of a sequence of i.i.d. random variables for the error we cannot use results as in Shao (1995) to obtain an almost sure bound for this process [cf. Tomkins (1974)]. But those results give an upper bound in probability (cf. Lemma A.2) for the maximum. This allows us to refine the bound above to $C \geq 2 + \delta$ for any $\delta > 0$ and obtain the rates of Theorem 6 in probability. We found that values of $C$ between 2 and 3 lead to good reconstruction for various simulation settings.

Figure 1 shows the behavior of the Potts minimizer for the test signals of Donoho and Johnstone (1994) sampled at 2048 points and a choice of $C = 2.5$. In order to understand the finite sample behavior of the Potts minimizer, the estimates are calculated at different signal-to-noise ratios $\|f\|^2/\sigma^2$ (seven, four and one). The reconstructions of the locally constant blocks signal (first row) differ very little from the original signal. This is not surprising since the original signal is in $S([0,1))$ where the estimator achieves parametric rates. The spikes of the bumps signal (second row) are correctly estimated for all cases. The estimator captures all relevant features of the Heavisine signal (third row) at the levels seven and four. Only in the presence of strong noise the detail of the spike right to the second maximum is lost. Finally, the case of the Doppler signal (fourth row) shows that the estimator adapts well to locally changing smoothness.

Clearly the performance depends on the particular function $f$. Hence one might want to try different approaches to selecting the parameter. One possibility is to choose the smoothing parameter according to the multiresolution criterion of Davies and Kovac (2001). If $f \in S([0,1))$, this criterion picks asymptotically the correct number of jumps.

THEOREM 5. *Assume $f \in S([0,1))$, $\xi_i^n \sim N(0, \sigma^2)$ i.i.d. and $\hat{\gamma}_n$ is chosen according to the MR-criterion, that is, $\hat{\gamma}_n$ is the maximal value such that the corresponding reconstruction $\hat{f}_n^{MR}$ satisfies*

$$(7) \qquad \frac{1}{\sqrt{\#I}} \left| \sum_{i \in I} Y_i^n - \hat{f}_n^{MR}(x_i^n) \right| \leq (1+\delta)\hat{\sigma}\sqrt{2 \log n}$$

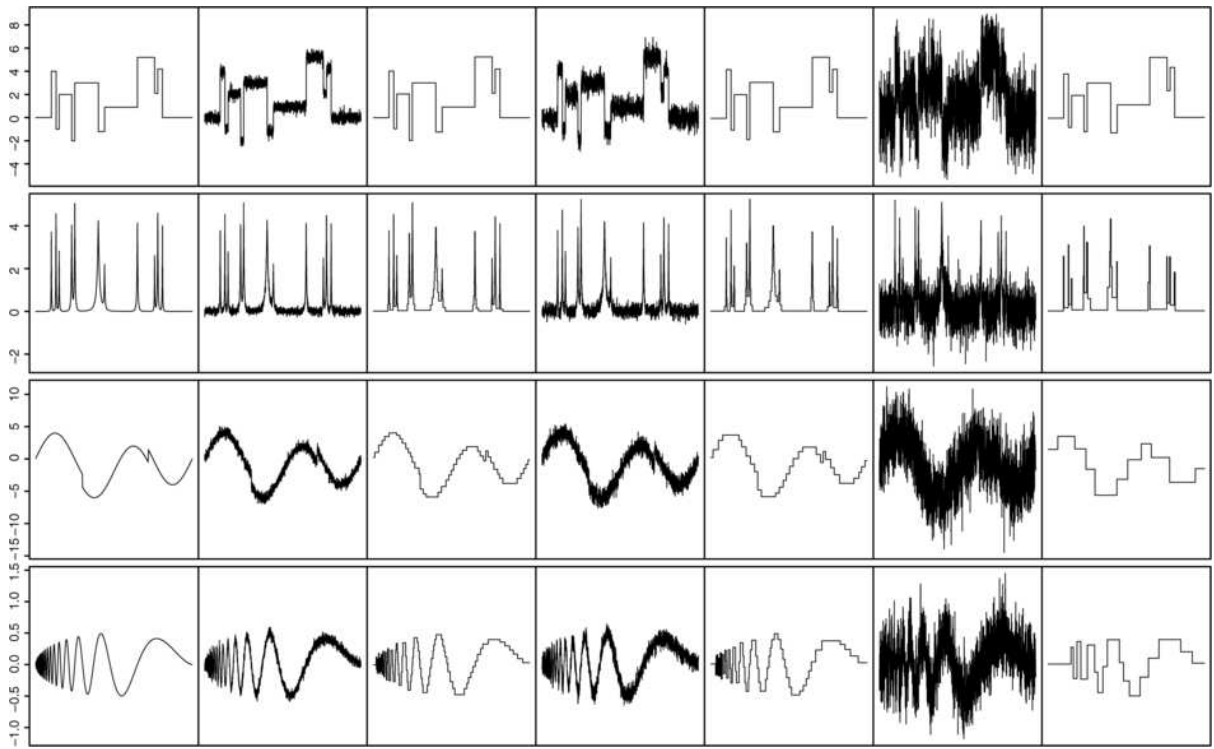

FIG. 1. *The left column shows signals from* Donoho and Johnstone (1994). *Columns 2, 4 and 6 show noisy versions with signal-to-noise ratios of 7, 4 and 1, respectively. On the right of each noisy signal is the Potts reconstruction. The penalty was chosen as* $\gamma_n = 2.5\hat{\sigma}^2 \log n/n$, *where* $\hat{\sigma}^2$ *is an estimate of the variance.*

*for all connected* $I \subset \{1, \ldots, n\}$, *some* $\delta > 0$ *and some consistent estimate* $\hat{\sigma}$ *of* $\sigma$. *Moreover, assume* $\gamma_n$ *satisfies condition* (H1) *and* $\hat{f}_n$ *is the corresponding reconstruction. Then* $\mathbb{P}(\hat{f}_n^{MR} = \hat{f}_n) \xrightarrow[n \to \infty]{} 1$.

Note that it is possible to derive the same result if in (7) only dyadic intervals [see Davies and Kovac (2001)] are considered. We conjecture that the MR-criterion leads to consistent estimates in more general settings.

**4. Discussion—relation to other models.** The Potts smoother falls in the general framework of van de Geer (2001) which gives very general and powerful tools to prove rates of convergence for penalized least squares estimates. With some effort, it is possible to use the methods developed in that paper to derive the convergence rates given in Theorem 2. However, using that method does not lead to the required constant in Section 3.4. In fact, the resulting constant in condition (H2) would be substantially larger.

Most penalized least squares methods either use a penalty which is a seminorm (as in spline regression) or penalizes the number or size of coefficients of an orthonormal basis reconstruction. Note that the Potts smoother belongs to none of these classes. Nonetheless, it is related to various other statistical procedures and we would like to close this paper by highlighting these relations and shortly comment on possible extensions to two dimensions.

*Bayesian interpretation and imaging.* In image analysis Bayesian methods for restoration have received much attention [see, e.g., Geman and Geman (1984)]. The Potts functional can be interpreted as a limit of the one-dimensional version of a certain MAP estimator, which has been used for edge-preserving smoothing, discussed by Blake and Zisserman (1987) and Künsch (1994) among many others. For a detailed discussion and overview of related functionals in dimension 1 [see Winkler et al. (2005)].

*Generalization to 2d.* For two-dimensional data, a measure of complexity corresponding to the number of jumps is given by the number of plateaus or partition elements. However, it is computationally infeasible to allow for arbitrary partitions in the reconstruction. Therefore one chooses a subclass of step functions with good approximation properties and seeks for effective minimization algorithms in this class. As in the one-dimensional case, the rate of convergence will be determined by the approximation properties of the chosen function class. One example, complexity penalized sums of squares with respect to a class of "Wedgelets" [cf. Donoho (1999)], is discussed in the Ph.D. thesis of Friedrich (2005), and possible alternatives in the survey by Führ, Demaret and Friedrich (2006). We mention that the proof of Theorem 2 could be adapted to their setting.

## APPENDIX: PROOFS

**A.1. Preliminaries.** Since the consistency results are formulated in terms of a function space, we translate all minimization problems to equivalent problems for functionals on $L^2([0,1))$. Therefore we introduce the functionals $\tilde{H}_\gamma^\infty(g,f) = H_\gamma^\infty(g,f) - \|f\|^2$ and $\tilde{H}_\gamma^n(g,f)$ is defined as $\tilde{H}_\gamma^\infty(g,f)$ for $g \in S_n([0,1)) := \iota^n(\mathbb{R}^n)$, and $\infty$, else. Clearly, the functionals are constructed in such a way that the minimization of $H_\gamma$ (3) on $\mathbb{R}^n$ is equivalent to the minimization of $\tilde{H}_\gamma^n$ if we identify the minimizers via the map $\iota^n$ defined in (4). The constant $-\|f\|^2$ is just added for convenience and does not affect the minimization. Obviously, $u \in \arg\min H_\gamma(\cdot, \overline{f}^n)$ if and only if $\iota^n(u) \in \arg\min \tilde{H}_\gamma^n(\cdot, f)$ and similarly for $H_\gamma(\cdot, y)$ for $y \in \mathbb{R}^n$. The most important property of these functionals is that the minimizers $g \in S([0,1))$ of $\tilde{H}_\gamma^n$ and $\tilde{H}_\gamma^\infty$ for $\gamma > 0$ are determined by their jump-set $J(g)$ and given by the projection onto the space of step functions which are constant outside that set. To make this precise in the course of the proofs, we introduce for any $J \subset (0,1)$ the partition $P_J = \{[a,b) : a, b \in J \cup \{0,1\}, (a,b) \cap J = \varnothing\}$. Abbreviating by

$$\mu_I(f) = \ell(I)^{-1} \int_I f(u)\, du$$

the mean of $f$ over some interval $I$, this projection is then given by

$$f_J = \sum_{I \in P_J} \mu_I(f) \mathbf{1}_I.$$

Further, we extend the noise in (1) to $L^2([0,1))$ by $\xi^n = \iota^n((\xi_1^n, \ldots, \xi_n^n))$ and, finally, we define for $f \in S([0,1))$ the minimum distance between any two jumps as

(8) $$\mathrm{mpl}(f) := \min\{|s-t| : s \neq t \in J(f) \cup \{0,1\}\}.$$

The proofs rely on properties of the noise, some a priori properties of the Potts minimizers and on proving epiconvergence of the functionals defined above with respect to the topology of $L^2([0,1))$.

**A.2. Two properties of the noise.** The behavior of $\xi_J^n = \sum_{I \in P_J} \mu_I(\xi^n)\mathbf{1}_I$ from Condition (A) is controlled by the following two estimates which are proved in Boysen et al. (2007), Section 4.2.

LEMMA A.1. *Let* $(\xi_i^n)_{n \in \mathbb{N}, 1 \leq i \leq n}$ *fulfill Condition (A). For*

(9) $$C_n := \sup_{1 \leq i \leq j \leq n} \frac{(\xi_i^n + \cdots + \xi_j^n)^2}{(j-i+1)\log n}$$

*we have that*

$$\limsup_{n\to\infty} C_n \leq 12\beta \qquad \mathbb{P}\text{-}a.s.$$

*Moreover, for all intervals $I \subset [0,1)$ and all $n \in \mathbb{N}$*

$$\mu_I(\xi^n)^2 \leq C_n \frac{\log n}{n\ell(I)}$$

*as well as*

$$(10) \qquad \|\xi_{J_n}^n\|^2 = \sum_{I \in P_{J_n}} \ell(I)\mu_I(\xi^n)^2 \leq C_n \frac{\log n}{n}(\#J_n + 1).$$

LEMMA A.2. *Assume $\xi_i^n \sim N(0,\sigma^2), i=1,\ldots,n$ i.i.d. for all $n$. Then for $C_n$ defined by (9) we have $C_n = 2\sigma^2 + o_\mathbb{P}(1)$.*

**A.3. A priori properties of the minimizers.** The following properties of the minimizers are used to prove our main statements.

LEMMA A.3. *Let $f \in L^2([0,1))$, $g \in \arg\min \tilde{H}_\gamma^n(\cdot, f)$ and $I \in P_{J(g)}$. Then, denoting $a = \mu_I(g) = \ell(I)^{-1}\int_I g(u)\,du$, the following statements are valid.*

(i) *If $I' \in P_{J(g)}$ and $I' \cup I$ is an interval, then*

$$\gamma \leq \frac{\ell(I)\ell(I')}{\ell(I) + \ell(I')}(\mu_I(f) - \mu_{I'}(f))^2.$$

(ii) *If $I' \in \mathcal{B}_n$, $I' \subset I$, is an interval, then*

$$2\gamma \geq \ell(I')(\mu_{I'}(f) - a)^2.$$

(iii) *If both $I' \in \mathcal{B}_n$ and $I' \cup I$ are intervals and $\mathbf{1}_{I'}g = b\mathbf{1}_{I'}$ for some $b \in \mathbb{R}$, then*

$$(b-a)\left(\mu_{I'}(f) - \frac{a+b}{2}\right) \geq 0.$$

(iv) *If $I_1', I_2', I_1' \cup I, I_2' \cup I \in \mathcal{B}_n$ are intervals and $\mathbf{1}_{I_l'}\hat{f} = b_l\mathbf{1}_{I_l'}$, $l=1,2$, then for all disjoint intervals $I_1, I_2 \in \mathcal{B}_n$, $I = I_1 \cup I_2$, such that $I_1 \cup I_1'$ and $I_2 \cup I_2'$ are intervals,*

$$\ell(I_1)(\mu_{I_1}(f) - b_1)^2 + \ell(I_2)(\mu_{I_2}(f) - b_2)^2$$
$$\geq \gamma + \ell(I_1)(\mu_{I_1}(f) - a)^2 + \ell(I_2)(\mu_{I_2}(f) - a)^2.$$

PROOF. The inequalities are obtained by elementary calculations comparing the values of $\tilde{H}_\gamma^n(\cdot, f)$ at $g$ and at some $\tilde{g}$ obtained from $g$ by: joining the plateaus at $I$ and $I'$ [for (i)], splitting the plateau at $I$ into three plateaus [for (ii)], moving the jump point [for (iii)], and removing the plateau at $I$ by joining each of the parts to the adjacent intervals [for (iv)].

As an example, we provide the calculations for (i). Determine $t$ by $\{t\} = \overline{I} \cap \overline{I'}$ and set $\tilde{g} = f_{J(g) \setminus \{t\}}$. Then $\tilde{g}$ differs from $g$ only on $I \cap I'$ such that

$$\begin{aligned}
0 &\leq \tilde{H}_\gamma^n(\tilde{g}, f) - \tilde{H}_\gamma^n(g, f) \\
&= -\gamma + \|(\mu_I(f) - \mu_{I \cup I'}(f))\mathbf{1}_I\|^2 + \|(\mu_{I'}(f) - \mu_{I \cup I'}(f))\mathbf{1}_{I'}\|^2 \\
&= -\gamma + \ell(I)(\mu_I(f) - \mu_{I \cup I'}(f))^2 + \ell(I')(\mu_{I'}(f) - \mu_{I \cup I'}(f))^2 \\
&= -\gamma + \frac{\ell(I)\ell(I')}{\ell(I) + \ell(I')}(\mu_I(f) - \mu_{I'}(f))^2,
\end{aligned}$$

which completes the proof of (i). $\square$

**A.4. Epiconvergence.** One basic idea of the consistency proofs is to use the concept of epiconvergence of the functionals [see, e.g., Dal Maso (1993), Hess (1996)]. We say that numerical functions $F_n : \Theta \mapsto \mathbb{R} \cup \{\infty\}$, $n = 1, \ldots, \infty$ on a metric space $(\Theta, \rho)$ epiconverge to $F_\infty$ if for all sequences $(\vartheta_n)_{n \in \mathbb{N}}$ with $\vartheta_n \to \vartheta \in \Theta$ we have $F_\infty(\vartheta) \leq \liminf_{n \to \infty} F_n(\vartheta_n)$, and for all $\vartheta \in \Theta$ there exists a sequence $(\vartheta_n)_{n \in \mathbb{N}}$ with $\vartheta_n \to \vartheta$ such that $F_\infty(\vartheta) \geq \limsup_{n \to \infty} F_n(\vartheta_n)$. One important property is that each accumulation point of a sequence of minimizers of $F_n$ is a minimizer of $F_\infty$. However, that does not mean that a sequence of minimizers has accumulation points at all. To prove this, one needs to show that the minimizers are contained in a compact set. The following lemma which is a straightforward consequence of the characterization of compact subsets of $D([0, 1))$ [Billingsley (1968), Theorem 14.3] will be applied to this end.

LEMMA A.4. *A subset $A \subset D([0, 1))$ is relatively compact if the following two conditions hold:*

(C1) *For all $t \in [0, 1]$ there is a compact set $K_t \subseteq \mathbb{R}$ such that*

$$g(t) \in K_t \quad \text{for all } g \in A.$$

(C2) *For all $\varepsilon > 0$ there exists a $\delta > 0$ such that for all $g \in A$ there is a step function $g_\varepsilon \in S([0, 1))$ such that*

$$\sup\{|g(t) - g_\varepsilon(t)| : t \in [0, 1]\} < \varepsilon \quad \text{and} \quad \mathrm{mpl}(g_\varepsilon) \geq \delta,$$

*where* mpl *is defined by (8).*

**A.5. The proof of Theorem 1(i), (ii) and Theorem 4(i), (ii).** For the sake of brevity we just give a short outline of the proof of the first two parts of Theorem 1 and the proof of Theorem 4(i). The details can be found in Boysen et al. (2007). The proof of Theorem 1(iii) is postponed to Section A.7, because it requires the proof of Theorem 3.

PROOF OF THEOREM 1(i), (ii). Note that condition (H1) automatically holds if $\gamma_n \to \gamma > 0$. We can thus prove both parts at once: Use first $\tilde{H}^n_{\gamma_n}(\hat{f}_n, f + \xi^n) \leq \tilde{H}^n_{\gamma_n}(0, f + \xi^n)$, $\gamma_n n / \log n \to \infty$ and (10) to obtain

$$\#J_n \leq \frac{2\|f\| + 2C_n(\log n/n)}{\gamma_n - 2C_n(\log n/n)} = O(\gamma_n^{-1}). \tag{11}$$

Then (10) and $\gamma_n n / \log n \to \infty$ imply

$$\|\xi^n_{J_n}\|^2 = \sum_{I \in P_{J_n}} \ell(I) \mu_I(\xi^n)^2 \to 0 \qquad \mathbb{P}\text{-a.s.} \tag{12}$$

The map

$$g \mapsto \begin{cases} \#J(g), & \text{if } g \in S([0,1)), \\ \infty, & \text{if } g \notin S([0,1)), \end{cases}$$

is lower semicontinuous as map from $L^2$ to $\mathbb{N} \cup \infty$. Using that together with (11) and (12), we can verify the two inequalities from the definition of epiconvergence and deduce that $\tilde{H}^n_{\gamma_n}(\cdot, f + \xi^n)$ actually converges to $\tilde{H}^\infty_\gamma(\cdot, f)$ for $\gamma_n \to \gamma \geq 0$ and $\gamma_n n / \log n \to \infty$ in that sense. Since for any $f \in L^2([0,1))$ the set $\{f_J : J \subset (0,1), \#J < \infty\}$ is relatively compact in $L^2([0,1))$, a comparison of $\tilde{H}^n_{\gamma_n}(\hat{f}_n, f + \xi^n)$ with $\tilde{H}^n_{\gamma_n}(0, f + \xi^n)$ and usage of (11) above yields that the set $\bigcup_{n \in \mathbb{N}} \arg\min \tilde{H}^n_{\gamma_n}(\cdot, f + \xi^n)$ is relatively compact. The uniqueness of the minimizer of $\tilde{H}^\infty_\gamma(\cdot, f)$ along with the epiconvergence of $\tilde{H}^n_{\gamma_n}(\cdot, f + \xi^n)$ and the compactness finally imply convergence of the minimizers. □

PROOF OF THEOREM 4(i). To prove this, one can proceed in a similar way as above. The proof of Lemma 1 is straightforward using $H^\infty_\gamma(0, f) = \|f\|^2$ and the relative compactness of $\{f_J : \#J \leq \|f\|^2/\gamma\}$ in $L^2([0,1))$ for $\gamma > 0$. □

Next, we will prove consistency in the space $D([0,1))$ equipped with the Skorokhod $J_1$-topology. This part is considerably more elaborate; in particular we need some of the a priori information about the minimizers provided by Lemma A.3.

PROOF OF THEOREM 4(ii). All equations in this proof hold $\mathbb{P}$-almost surely, which will be omitted for ease of notation. If $f_1, f_2 \in D([0,1))$ are limit points of the sequence of minimizers, we know by Theorem 1(ii) that $f = f_1 = f_2$ in $L^2([0,1))$, which implies that they are equal in $D([0,1))$. Thus, it is enough to show that the minimizers $\{(f + \xi^n)_{J_n} : n \in \mathbb{N}\}$ are contained in a compact set. For this goal we use now the conditions (C1), (C2) from Lemma A.4.

For the proof of (C1), consider any interval $I \in P_{J_n}$. We know from part (i) of Lemma A.3, for any neighboring interval $I'$, that

$$\gamma_n \leq \frac{\ell(I)\ell(I')}{\ell(I) + \ell(I')}(\mu_I(f + \xi^n) - \mu_{I'}(f + \xi^n))^2$$

$$\leq \frac{\ell(I)\ell(I')}{\ell(I) + \ell(I')}\left(12\|f\|_\infty^2 + 3C_n\frac{\log n}{n\ell(I)} + 3C_n\frac{\log n}{n\ell(I')}\right)$$

$$\leq 12\|f\|_\infty^2 \ell(I) + 6C_n\frac{\log n}{n}.$$

This yields $1/\ell(I) = O(\gamma_n^{-1})$. Application of Lemma A.1 yields

$$\|\xi^n_{J_n}\|_\infty^2 = \max\{\mu_I(\xi^n)^2 : I \in P_{J_n}\} = O\left(\frac{\log n}{n\gamma_n}\right) = o(1)$$

and $\|(f + \xi^n)_{J_n}\|_\infty = O(1)$. For the proof of (C2), let us fix $\varepsilon > 0$ and a step function $\tilde{f}$ with $\|f - \tilde{f}\|_\infty < \varepsilon/7$. Further, set $\delta = \mathrm{mpl}(\tilde{f}) > 0$. Now we will consider three different classes of intervals $I \in P_{J_n}$ which are characterized by their position relative to $J(\tilde{f})$ and estimate $(f + \xi^n)_{J_n} - \tilde{f}$ uniformly on them, separately.

*Class* 1 consists of intervals $I$ with $J(\tilde{f}) \cap I = \varnothing$. We obtain that

$$\|\mathbf{1}_I(\tilde{f} - (f + \xi^n)_{J_n})\|_\infty \leq \|\mathbf{1}_I(\tilde{f} - f_{J_n})\|_\infty + \|\xi^n_{J_n}\|_\infty \leq \|\tilde{f} - f\|_\infty + o(1) < \varepsilon/7$$

for large enough $n$ uniformly for all such $I$ and $n$.

*Class* 2 covers intervals $I$ which are not in class 1 but for which there is some interval $\tilde{I} \in P_{J(\tilde{f})}$ with $\ell(I \cap \tilde{I}) \geq \delta/6$. To apply Lemma A.3(ii), choose an interval $I' \subseteq I \cap \tilde{I}$ from $\mathcal{B}_n$ such that $\rho_H(I', I \cap \tilde{I}) \leq 1/n$. We find for all $t \in I'$

$$|(f + \xi^n)_{J_n}(t) - \mu_{I'}(f + \xi^n)| \leq \sqrt{\frac{2\gamma_n}{\ell(I')}} \leq \sqrt{\frac{2\gamma_n}{\delta/6 - 2/n}},$$

hence

$$|(f + \xi^n)_{J_n}(t) - \tilde{f}(t)| \leq |\mu_{I'}(f) - \mu_{I'}(\tilde{f})| + |\mu_{I'}(\xi^n)| + \sqrt{\frac{2\gamma_n}{\delta/6 - 2/n}}$$

$$\leq \varepsilon/7 + \sqrt{\frac{C_n \log n/n}{\delta/6 - 2/n}} + \sqrt{\frac{2\gamma_n}{\delta/6 - 2/n}} < \varepsilon/6$$

for large enough $n$ depending only on $(\gamma_n)_{n\in\mathbb{N}}, \delta, \varepsilon$. Clearly, this implies that for $n$ large enough $\sup_{I\cap I'} |(f+\xi^n)_{J_n} - \tilde{f}| < \varepsilon/6$ uniformly in $I, I'$.

*Class* 3 contains all intervals $I \in P_{J_n}$ which are in neither class 1 nor class 2 such that $\ell(I) < \delta/3$ and $I \cap J(\tilde{f}) = \{t_0\}$. Then the neighboring intervals of $I$ in $P_{J_n}$ belong necessarily to class 1 or 2. Further, if a neighboring interval $I'$ is in class 2, we know that there is $\tilde{I} \in P_{J(\tilde{f})}$ with $\ell(\tilde{I} \cap I') \geq \delta/6$ and $\tilde{I} \cap I \neq \varnothing$ such that $\text{dist}(t_0, \tilde{I}) = 0$. In any case, we find for any interval $\tilde{I}$ with endpoint $t_0$ in $P_{J(\tilde{f})}$ and any interval $I'$ neighboring $I$ in $P_{J_n}$ with $I' \cap \tilde{I} \neq \varnothing$ that $\sup_{\tilde{I}\cap I'} |(f+\xi^n)_{J_n} - \tilde{f}| < \varepsilon/6$ and thus $|\mu_{I'}((f+\xi^n)_{J_n}) - \mu_{\tilde{I}}(\tilde{f})| = |\mu_{\tilde{I}\cap I'}((f+\xi^n)_{J_n}) - \mu_{\tilde{I}\cap I'}(\tilde{f})| < \varepsilon/6$.

We choose $t_1$ with $nt_1 \in \mathbb{N}$ and $|t_1 - t_0| < 1/n$ as well as $I_1 = I \cap [0, t_1)$, $I_2 = I \cap [t_1, 1)$ and $I'_j$ as neighboring intervals of $I_j$ in $P_{J_n}$, $j = 1, 2$. Denoting $a = \mu_I(f+\xi^n)$ and $b_j = \mu_{I'_j}(f+\xi^n)$, application of Lemma A.3(iv) yields (together with Lemma A.1) that

$$\ell(I_1)(a - \mu_{I_1}(f+\xi^n))^2 + \ell(I_2)(a - \mu_{I_2}(f+\xi^n))^2$$
$$\leq -\gamma_n + \ell(I_1)(b_1 - \mu_{I_1}(f+\xi^n))^2 + \ell(I_2)(b_2 - \mu_{I_2}(f+\xi^n))^2,$$
$$\ell(I_1)(a - \mu_{I_1}(f))^2 + \ell(I_2)(a - \mu_{I_2}(f))^2$$
$$\leq -\gamma_n + 2\ell(I_1)\mu_{I_1}(\xi^n)(a-b_1) + 2\ell(I_2)\mu_{I_2}(\xi^n)(a-b_2)$$
$$\quad + \ell(I_1)(b_1 - \mu_{I_1}(f))^2 + \ell(I_2)(b_2 - \mu_{I_2}(f))^2$$
$$\leq 2\ell(I_1)\mu_{I_1}(\xi^n)(a-b_1) + 2\ell(I_2)\mu_{I_2}(\xi^n)(a-b_2) + \ell(I)\varepsilon^2(1/6 + 1/7)^2$$
$$\leq 2|a - b_1|\sqrt{\ell(I_1)C_n \log n/n} + 2|a - b_2|\sqrt{\ell(I_2)C_n \log n/n} + \ell(I)\varepsilon^2/9.$$

From $\|\xi^n_{J_n}\| = o(1)$ we find $b_i - a = O(1)$ such that for large $n$ depending on $\varepsilon, \delta$ only

$$\ell(I_1)(a - \mu_{I_1}(f))^2 + \ell(I_2)(a - \mu_{I_2}(f))^2 \leq \ell(I)\varepsilon^2/9.$$

The above results yield for $t' \in I$ that

$$\ell(I_1)((f+\xi^n)_{J_n}(t') - \mu_{I_1}(f))^2 + \ell(I_2)((f+\xi^n)_{J_n}(t') - \mu_{I_2}(f))^2 \leq \ell(I)\varepsilon^2/9$$

and hence

$$\min(|(f+\xi^n)_{J_n}(t') - \mu_{I_1}(f)|, |(f+\xi^n)_{J_n}(t') - \mu_{I_2}(f)|) \leq \varepsilon/3,$$
$$\min(|(f+\xi^n)_{J_n}(t') - \mu_{I_1}(\tilde{f})|, |(f+\xi^n)_{J_n}(t') - \mu_{I_2}(\tilde{f})|) \leq \varepsilon/2.$$

This shows that either $\|\mathbf{1}_{I\cap[t_0,1)}(\tilde{f} - (f+\xi^n)_{J_n})\|_\infty \leq \varepsilon/2$ or $\|\mathbf{1}_{I\cap[0,t_0)}(\tilde{f} - (f+\xi^n)_{J_n})\|_\infty \leq \varepsilon/2$ holds for large $n$, depending on $\varepsilon, \delta$ only.

Given $J_n$ we define a new partition $P'_n$ coarser than $P_{J_n}$ by the following procedure. First we join all neighboring intervals of class 1 and denote the resulting intervals again as class 1. If there are class 1 intervals left of length $< \delta/3$, there must be a left or a right neighbor which is class 2 and has an overlap of length $> \delta/3$ with an interval of constancy of $\tilde{f}$. Then we join the class 1 interval to that neighbor (if there are two, to the left one). At the end, we join each class 3 interval $I$ to its left neighbor, if $\|\mathbf{1}_{I \cap [t_0,1)}(\tilde{f} - (f + \xi^n)_{J_n})\|_\infty \leq \varepsilon/2$, or else to its right neighbor. The collection of those joined intervals is $P'_n$.

By the results for class 1, 2, 3 intervals we know for all $I \in P'_n$ that $\ell(I) \geq \delta/3$. Further, for each $I \in P'_n$ there is $I' \in P_{J(\tilde{f})}$ such that $\tilde{I} \cap I' \neq \varnothing$ for all $\tilde{I} \in P_{J_n}$, $\tilde{I} \subseteq I$, and $\|\mathbf{1}_{I \cap I'}(\tilde{f} - (f + \xi^n)_{J_n})\|_\infty < \varepsilon/2$ holds. Thus, defining $\tilde{f}_n = \sum_{I \in P'_n} \mu_I((f + \xi^n)_{J_n})\mathbf{1}_I$ we obtain that $\|\tilde{f}_n - (f + \xi^n)_{J_n}\|_\infty < \varepsilon$. Thus (C2) is established and by Lemma A.4 $\{(f + \xi^n)_{J_n} : n \in \mathbb{N}\}$ is contained in a compact set. This completes the proof of the first assertion. The second assertion follows from the fact that convergence in $D([0, 1))$ implies convergence in $L_\infty([0, 1])$ if the limit is continuous [Billingsley (1968), page 112]. □

**A.6. The proof of Theorem 2.** Fix numbers $k_n \geq 1$, the precise magnitude of which will be chosen below. Further, sets $K_n \subseteq \{1/n, \ldots, (n-1)/n\}$ are chosen such that $f_{K_n}$ is a best approximation of $f$ by a step function from $S_n([0, 1))$ with $k_n \geq 1$ jumps, which exists since the subspace of $S_n([0, 1))$ containing functions $g$ with $\#J(g) \leq k_n$ and $\|g\| \leq 2\|f\|$ is compact.

Let $\tilde{f}_{k_n}$ be an approximation of $f$ in $S([0, 1))$ with at most $k_n$ jumps for which $\|\tilde{f}_{k_n} - f\| = O(\frac{1}{k_n^\alpha})$. Further, without loss of generality, we can assume that $\tilde{f}_{k_n} = f_{J(\tilde{f}_{k_n})}$ which implies $\|\tilde{f}_{k_n}\|_\infty \leq \|f\|_\infty$. Moving each jump of $\tilde{f}_{k_n}$ to the next $t \in [0, 1]$ with $nt \in \mathbb{N}$ but leaving the value of $\tilde{f}_{k_n}$ unchanged on each plateau, we obtain a step function $\tilde{f}_n \in S_n([0, 1))$ with $\|\tilde{f}_{k_n} - \tilde{f}_n\|^2 \leq \frac{2k_n}{n}\|f\|_\infty^2$. This shows $\|\tilde{f}_n - f\|^2 = O(\frac{1}{k_n^{2\alpha}} + \frac{k_n}{n})$. Since $f_{K_n}$ is a best approximation, we derive

$$\|f_{K_n} - f\|^2 = O\left(\frac{1}{k_n^{2\alpha}} + \frac{k_n}{n}\right).$$

By definition $\hat{f}_n$ is a minimizer of $\tilde{H}^n_{\gamma_n}(\cdot, f + \xi^n)$ and we get

$$\tilde{H}^n_{\gamma_n}(\hat{f}_n, f + \xi^n) \leq \tilde{H}^n_{\gamma_n}(f_{K_n}, f + \xi^n).$$

By $\#K_n = k_n$, this implies $\gamma_n \#J_n + \|\hat{f}_n - f - \xi^n\|^2 \leq \gamma_n k_n + \|f_{K_n} - f - \xi^n\|^2$ and hence

$$\|\hat{f}_n - f\|^2 \leq \gamma_n(k_n - \#J_n) + \|f_{K_n} - f\|^2 + 2\langle f - f_{K_n}, \xi^n\rangle + 2\langle \hat{f}_n - f, \xi^n\rangle$$
$$\leq \gamma_n(k_n - \#J_n) + \|f_{K_n} - f\|^2 + 2\langle \hat{f}_n - f_{K_n}, \xi^n\rangle.$$

Now observe that $J(\hat{f}_n - f_{K_n}) \subseteq J_n \cup K_n$ which gives

$$\begin{aligned}\langle \hat{f}_n - f_{K_n}, \xi^n \rangle &= \langle \hat{f}_n - f_{K_n}, (\xi^n)_{J_n \cup K_n} \rangle \\ &\leq \|\hat{f}_n - f_{K_n}\| \|(\xi^n)_{J_n \cup K_n}\| \\ &\leq \|\hat{f}_n - f\| \|(\xi^n)_{J_n \cup K_n}\| + \|f - f_{K_n}\| \|(\xi^n)_{J_n \cup K_n}\| \\ &\leq \frac{1}{2+\delta}\|\hat{f}_n - f\|^2 + \frac{2+\delta}{4}\|(\xi^n)_{J_n \cup K_n}\|^2 \\ &\quad + \frac{1}{\delta}\|f - f_{K_n}\|^2 + \frac{\delta}{4}\|(\xi^n)_{J_n \cup K_n}\|^2.\end{aligned}$$

The above inequalities yield

$$\frac{\delta}{2+\delta}\|\hat{f}_n - f\|^2 \leq \gamma_n(k_n - \#J_n) + \frac{2+2\delta}{\delta}\|f_{K_n} - f\|^2 + (1+\delta)\|(\xi^n)_{J_n \cup K_n}\|^2.$$

Using the estimate (10) with $C_n$ from (9) we obtain for $C' = \delta/(2+\delta)$

$$\begin{aligned}C'\|\hat{f}_n - f\|^2 &\leq \gamma_n(k_n - \#J_n) + C''\left(\frac{1}{k_n^{2\alpha}} + \frac{k_n}{n}\right) + (1+\delta)C_n \frac{\log n}{n}(\#J_n + k_n + 1) \\ &\leq k_n\left(\gamma_n + (1+\delta)C_n \frac{\log n}{n} + \frac{C''}{n}\right) + \#J_n\left((1+\delta)C_n \frac{\log n}{n} - \gamma_n\right) \\ &\quad + \frac{C''}{k_n^{2\alpha}} + (1+\delta)C_n \frac{\log n}{n},\end{aligned}$$

for some constant $C''$ depending on $f$. We get from $\gamma_n \geq (1+\delta)12\beta \log n/n$ together with the relation $\limsup_{n\to\infty} C_n \leq 12\beta$ that $(1+\delta)C_n \log n/n \leq \gamma_n$ and $C''/n \leq \gamma_n$ for large enough $n$, hence $C'\|\hat{f}_n - f\|^2 \leq \gamma_n(3k_n + 1) + C''/k_n^{2\alpha}$. Choosing $k_n = \lfloor \gamma_n^{-1/(2\alpha+1)} \rfloor$ we obtain

$$\|\hat{f}_n - f\|^2 = O(\gamma_n^{2\alpha/(2\alpha+1)})$$

and the proof is complete.

### A.7. The proof of Theorem 3, Theorem 1(iii) and Theorem 4(iii).

PROOF OF THEOREM 3(ii). 1. First we will show that

(13) $\quad \forall t \in J(f) \; \exists t_n \in J_n \quad$ with $|t_n - t| < \mathrm{mpl}(f)/3.$

From part (i) of Theorem 4 and $S([0,1)) \subset D([0,1))$ we obtain immediately that $\hat{f}_n \xrightarrow[n\to\infty]{D([0,1))} f$. Therefore, there is some random integer $n_0$ such that for

all $n \geq n_0$

(14)
$$\rho_S(\hat{f}_n, f) < \min(\min\{|f(t) - f(t-0)| : t \in J(f)\}/2, |\log(1 - \tfrac{2}{3}\mathrm{mpl}(f))|).$$

The relation (13) is a direct consequence of inequality (14). Assume (13) does not hold. In this case, a Lipschitz function $\lambda \in \Lambda_1$ with $L(\lambda) < |\log(1 - 2/3\mathrm{mpl}(f))|$ could not achieve $t \in J(\hat{f}_n \circ \lambda)$ and hence $\|\hat{f}_n \circ \lambda - f\|_\infty \geq |f(t) - f(t-0)|/2$ contradicting (14).

2. Now we will show that for all $t \in J(f)$ there exists a sequence $t_n \in J_n$, such that $|t_n - t| = O(\log n/n)$. For any $t \in J(f)$ let $t_n$ be a point in $J_n$ closest to $t$. We want to apply Lemma A.3(iii). For that goal, suppose for the moment that $t_n < t$ and $f(t) > f(t-0)$. Choose $I_n \in P_{J_n}$ as interval with right end point $t_n$ and set $I'_n = [t_n, s_n)$ where $ns_n \in \mathbb{N}$ is such that $|s_n - t| < 1/n$ as well as $a_n = \mu_{I_n}(\hat{f}_n)$ and $b_n = \mu_{I'_n}(\hat{f}_n)$. Then Lemma A.3(iii) shows

$$(b_n - a_n)\left(\mu_{I'_n}(f + \xi^n) - \frac{a_n + b_n}{2}\right) \geq 0.$$

Clearly, $\hat{f}_n \xrightarrow[n \to \infty]{D([0,1))} f$ implies $a_n \xrightarrow[n \to \infty]{} f(t-0)$ and $b_n \xrightarrow[n \to \infty]{} f(t)$ such that almost surely eventually

$$\mu_{I'_n}(f) - \frac{a_n + b_n}{2} \geq -\mu_{I'_n}(\xi^n) \geq -C_n\sqrt{\frac{\log n}{n\ell(I'_n)}}.$$

We know further $\lim_{n \to \infty} \mu_{I'_n}(f) = f(t-0)$ such that almost surely eventually

$$0 > \frac{f(t-0) - f(t)}{3} \geq -C_n\sqrt{\frac{\log n}{n\ell(I'_n)}}$$

which implies $\ell(I'_n) = O(\log n/n)$ and $|t_n - t| = O(\log n/n)$.

3. Next we will prove that there exists no sequence $t_n \in J_n$ which satisfies the relation $\limsup_{n \to \infty}(n/\log n)\rho_H(\{t_n\}, J) = \infty$. We consider two adjacent intervals $I, I' \in P_{J_n}$ for which there is an $\tilde{I} \in P_{J(f)}$ with $\ell(I \cup I' \setminus \tilde{I}) = O(\log n/n)$. Then

$$|\mu_I(f) - \mu_{I \cap \tilde{I}}(f)| = \frac{|\ell(I \cap \tilde{I})\int_I f(u)\,du - \ell(I)\int_{I \cap \tilde{I}} f(u)\,du|}{\ell(I)\ell(I \cap \tilde{I})}$$

$$= \frac{|\ell(I \cap \tilde{I})\int_{I \setminus \tilde{I}} f(u)\,du - \ell(I \setminus \tilde{I})\int_{I \cap \tilde{I}} f(u)\,du|}{\ell(I)\ell(I \cap \tilde{I})}$$

$$\leq 2\|f\|_\infty \frac{\ell(I \cap \tilde{I})\ell(I \setminus \tilde{I})}{\ell(I)\ell(I \cap \tilde{I})} = 2\|f\|_\infty \frac{\ell(I \setminus \tilde{I})}{\ell(I)}$$

and a similar estimate holds for $I'$. By means of $\mu_{I \cap \tilde{I}}(f) = \mu_{I' \cap \tilde{I}}(f)$ and

$1/\ell(I) = O(1/\gamma_n)$ we obtain

$$(\mu_I(f) - \mu_{I'}(f))^2 \leq (1/\ell(I)^2 + 1/\ell(I')^2)O\left(\frac{\log^2 n}{n^2}\right)$$

$$\leq (1/\ell(I) + 1/\ell(I'))O\left(\frac{\log^2 n}{\gamma_n n^2}\right)$$

$$= (1/\ell(I) + 1/\ell(I'))o\left(\frac{\log n}{n}\right).$$

Now Lemma A.3(i) implies

$$\gamma_n \leq \frac{\ell(I)\ell(I')}{\ell(I) + \ell(I')}(\mu_I(f + \xi^n) - \mu_{I'}(f + \xi^n))^2$$

$$\leq \frac{\ell(I)\ell(I')}{\ell(I) + \ell(I')}(3(\mu_I(f) - \mu_{I'}(f))^2 + 3\mu_I(\xi^n)^2 + 3\mu_{I'}(\xi^n)^2)$$

$$\leq O\left(\frac{\log n}{n}\right)\frac{\ell(I)\ell(I')}{\ell(I) + \ell(I')}(1/\ell(I) + 1/\ell(I')) = O\left(\frac{\log n}{n}\right).$$

This contradicts $\gamma_n n / \log n \to \infty$. Thus, almost surely, there are only finitely many $n$ for which there are two adjacent intervals $I, I' \in P_{J_n}$ and $\tilde{I} \in P_{J(f)}$ with $\ell(I \cup I' \setminus \tilde{I}) = O(\log n/n)$. Consequently, $\rho_H(J_n, J(f)) = O(\log n/n)$, which implies the statement. $\square$

PROOF OF THEOREM 3(i). 4. Suppose now there are $s_n, t_n \in J_n$ with $s_n \to t$, $t_n \to t$ for $t \in J(f)$. Then we have by the previous result that $|t_n - s_n| = O(\log n/n)$ as well as $1/|t_n - s_n| = O(1/\gamma_n)$. This gives us $\log n/(n\gamma_n) = O(1)$ contradicting $n\gamma_n/\log n \to \infty$. Thus $\#J_n = \#J(f)$ eventually. $\square$

PROOF OF THEOREM 3(iii). 5. For this statement, observe that in the special situation considered in step 2, it is not necessary to assume $|s_n - t| < 1/n$. Hence for any $s_n \in [t_n, t)$ with $ns_n \in \mathbb{N}$ we have almost surely eventually

$$0 > \frac{f(t-0) - f(t)}{3} \geq -\mu_{[t_n, s_n)}(\xi^n)$$

conditional on $t_n < t$. Denote $p$ the largest integer such that $p/n \leq t - 1/n$. Using the exponential inequality [cf. Petrov (1975), Sections 3 and 4]

$$(15) \qquad \mathbb{P}\left(\sum_{i=1}^{n} \mu_i \xi_i^n \geq z\right) \leq \exp\left(-\frac{z^2}{4\beta \sum_{i=1}^{n} \mu_i^2}\right) \qquad \text{for all } z \in \mathbb{R},$$

for triangular arrays fulfilling Condition (A) and all numbers $\mu_i$, $i = 1, \ldots, n$,

we obtain for all $k' \in \mathbb{N}$

$$\mathbb{P}(\{k'/n < (t - t_n) \leq (k' + 1)/n\})$$
$$\leq \mathbb{P}\bigg(\bigg\{\mu_{[(p+1-k')/n,(p+1-k'+i)/n)}(\xi^n) \geq \frac{f(t) - f(t-0)}{3}$$
$$\text{for all } i = 1, \ldots, k'\bigg\}\bigg)$$
$$= \mathbb{P}\bigg(\bigg\{\frac{\xi^n_{p-k'+1} + \cdots + \xi^n_{p-k'+i}}{i} \geq \frac{f(t) - f(t-0)}{3} \text{ for all } i = 1, \ldots, k'\bigg\}\bigg)$$
$$\leq \mathbb{P}\bigg(\bigg\{\frac{\xi^n_{p-k'+1} + \cdots + \xi^n_p}{k'} \geq \frac{f(t) - f(t-0)}{3}\bigg\}\bigg)$$
$$\leq \exp\bigg(\frac{-k'z^2}{4\beta}\bigg) = \bigg(\exp\bigg(\frac{-z^2}{4\beta}\bigg)\bigg)^{k'} =: q^{k'},$$

where $z = (f(t) - f(t-0))/3$. Note that $q < 1$ depends on $f(t) - f(t-0)$ and $\beta$ only. Clearly, we can use a similar argument if $f(t-0) > f(t)$ or $t_n \geq t$. Summing up these inequalities we obtain $\mathbb{P}(\{|t - t_n| \geq k/n\}) \leq 2q^k/(1-q)$ and

$$\mathbb{P}(\{\rho_H(J_n, J(f)) \geq k/n\}) \leq 2\#J(f)q^k/(1-q).$$

This shows $\lim_{k \to \infty} \limsup_{n \to \infty} \mathbb{P}(\{\rho_H(J_n, J(f)) \geq k/n\}) = 0$, or in other words $\rho_H(J_n, J(f)) = O_\mathbb{P}(n^{-1})$. $\square$

PROOF OF THEOREM 1(iii), THEOREM 4(iii). 6. By 4 and 5, we may choose $n$ so large that $\#J_n = \#J(f)$ and $\rho_H(J_n, J(f)) \leq \text{mpl}(f)/3$. Then there is a unique 1–1 map $\varphi_n: J(f) \longmapsto J_n$ for which $\sum_{t \in J(f)} |t - \varphi_n(t)|$ is minimal. We derive $\varphi_n(t) - t = O(\log n/n)$ for all $t \in J(f)$. Extend now $\varphi_n$ by $\varphi_n(0) = 0$ and $\varphi_n(1) = 1$. For $[s, t] \in P_{J(f)}$ we get thus

$$\|\mathbf{1}_{[\varphi_n(s),\varphi_n(t))} - \mathbf{1}_{[s,t)}\| = O\bigg(\sqrt{\frac{\log n}{n}}\bigg).$$

Further, $\|f\|_\infty < \infty$ yields $|\mu_{[\varphi_n(s),\varphi_n(t))}(f) - \mu_{[s,t)}(f)| = O(\sqrt{\log n/n})$. Lemma A.1 implies that $\mu_{[\varphi_n(s),\varphi_n(t))}(\xi^n) = O(\sqrt{\log n/n})$ such that $\|\hat{f}_n - f\| = O(\sqrt{\log n/n})$ which yields the first part of Theorem 1(iii) and

$$\|\mu_{[\varphi_n(s),\varphi_n(t))}(f + \xi^n)\mathbf{1}_{[\varphi_n(s),\varphi_n(t))} - \mu_{[s,t)}(f)\mathbf{1}_{[s,t)}\| = O\bigg(\sqrt{\frac{\log n}{n}}\bigg).$$

We define an extension $\lambda_n \in \Lambda_1$ of $\varphi_n$ by linear interpolation. From above, we obtain the estimate $\|\hat{f}_n - f \circ \lambda_n\|_\infty = O(\sqrt{\log n/n})$. Furthermore,

$$\text{L}(\varphi_n) = \max_{[s,t] \in P_{J(f)}} \bigg|\log \frac{\varphi_n(t) - \varphi_n(s)}{t - s}\bigg| = O(\log n/n)$$

such that $\rho_S(\hat{f}_n, f) = O(\sqrt{\log n/n})$.

7. By direct calculations we obtain from (15) and Lemma A.1 that

$$\max_{I \in P_{J_n}} |\mu_I(\xi^n)| = O_\mathbb{P}(n^{-1/2}).$$

Using this estimate and $\rho_H(J_n, J(f)) = O_\mathbb{P}(\frac{1}{n})$ in the same way as the almost sure rate in step 6, we obtain that $\rho_S(\hat{f}_n, f)$ and $\|\hat{f}_n - f\|$ are of order $O_\mathbb{P}(1/\sqrt{n})$. □

**A.8. The proof of Theorem 5.** It is sufficient to show that

$$\mathbb{P}(\#J(\hat{f}_n^{MR}) = \#J(\hat{f}_n)) \xrightarrow[n \to \infty]{} 1.$$

Assume there exists some subsequence $n_k$ such that $\#J(\hat{f}_{n_k}^{MR}) < \#J(f)$ for all $n_k$. As a step function with $\#J(f)$ jumps cannot be approximated by a sequence of functions with fewer jumps, there exists a sequence of connected intervals $I_{n_k}$ with $I_{n_k} \in \mathcal{B}_{n_k}$ such that $\liminf_{n_k \to \infty} l(I_{n_k}) \geq \epsilon_1 > 0$ and for $\tilde{I}_{n_k} = \{i : x_i^{n_k} \in I_{n_k}\}$

$$\left| \frac{1}{\#\tilde{I}_{n_k}} \sum_{i \in \tilde{I}_{n_k}} \overline{f}_i^{n_k} - \hat{f}_{n_k}^{MR}(x_i^{n_k}) \right| \geq \epsilon_2 > 0.$$

Consequently by Lemma A.1 for large $n_k$

$$\frac{|\sum_{i \in \tilde{I}_{n_k}} Y_i^{n_k} - \hat{f}_{n_k}^{MR}(x_i^{n_k})|}{\sqrt{\#\tilde{I}_{n_k}}}$$

$$\geq \epsilon_2 \sqrt{\epsilon_1 n_k} - \frac{|\sum_{i \in \tilde{I}_{n_k}} \xi_i^{n_k}|}{\sqrt{\#\tilde{I}_{n_k}}}$$

$$\geq \epsilon_2 \sqrt{\epsilon_1 n_k} - O(\sqrt{\log n_k}) \qquad \mathbb{P}\text{-a.s.}$$

This implies that for large $n_k$ the MR-criterion is not satisfied. By Theorem 3(i) we have $\mathbb{P}(\#J(\hat{f}_n) = \#J(\hat{f})) \to 1$ for $n \to \infty$. Hence $\mathbb{P}(\#J(\hat{f}_n^{MR}) \geq \#J(\hat{f}_n)) \xrightarrow[n \to \infty]{} 1$.

It remains to show that $\hat{f}_n^{MR}$ has asymptotically at most as many jumps as $\hat{f}_n$. Observe that

(16)
$$\max_{1 \leq j \leq k \leq n} \frac{|\sum_{i=j}^k Y_i^n - \hat{f}_n(x_i^n)|}{\sqrt{k-j+1}}$$
$$\leq \max_{1 \leq j \leq k \leq n} \frac{|\sum_{i=j}^k \xi_i^n| + |\sum_{i=j}^k \hat{f}_n(x_i^n) - \overline{f}_i^n|}{\sqrt{k-j+1}}.$$

By the Cauchy–Schwarz inequality and Theorem 1(iii) we have for $1 \leq j \leq k \leq n$

$$\frac{|\sum_{i=j}^k \hat{f}_n(x_i^n) - \overline{f}_i^n|}{\sqrt{k-j+1}} = \frac{n|\langle \hat{f}_n - \overline{f}^n, \mathbf{1}_{[j/n,(k+1)/n)} \rangle|}{\sqrt{k-j+1}}$$

$$\leq n \frac{\|\mathbf{1}_{[j/n,(k+1)/n)}\|}{\sqrt{k-j+1}} (\|\hat{f}_n - f\| + \|\overline{f}^n - f\|)$$

$$= \sqrt{n}(\|\hat{f}_n - f\| + \|\overline{f}^n - f\|) = O_P(1)$$

uniformly in $j, k$. Lemma A.2 implies

$$\max_{1 \leq j \leq k \leq n} \frac{|\sum_{i=j}^k \xi_i^n|}{\sqrt{k-j+1}} = \sigma\sqrt{2\log n} + o_P(\sqrt{\log n}).$$

Applying the results above to (16) we arrive at

$$\max_{1 \leq j \leq k \leq n} \frac{|\sum_{i=j}^k Y_i^n - \hat{f}_n(x_i^n)|}{\sqrt{k-j+1}} = \sigma\sqrt{2\log n} + o_P(\sqrt{\log n}).$$

Since $\hat{\sigma}$ is a consistent estimate of $\sigma$, this implies that the probability that $\hat{f}_n$ satisfies the MR-criterion tends to 1 as $n$ goes to infinity. As $\hat{\gamma}_n$ is chosen maximal such that the MR-criterion is satisfied, we can conclude $\mathbb{P}(\hat{\gamma}_n \geq \gamma_n) \xrightarrow[n \to \infty]{} 1$ and consequently $\mathbb{P}(\#J(\hat{f}_n^{MR}) \leq \#J(\hat{f}_n)) \xrightarrow[n \to \infty]{} 1$ which proves the claim. $\square$

**Acknowledgment.** We wish to thank L. Birgé, L. Brown, L. Dümbgen, F. Friedrich, K. Gröchenig, T. Hotz, E. Liebscher, E. Mammen, G. Winkler, two referees and two Associate Editors for helpful comments and bibliographic information.

## REFERENCES

Aurich, V. and Weule, J. (1995). Nonlinear Gaussian filters performing edge preserving diffusion. In *Proc. 17. DAGM-Symposium, Bielefeld* 538–545. Springer, Berlin.

Billingsley, P. (1968). *Convergence of Probability Measures*. Wiley, New York. MR0233396

Birgé, L. and Massart, P. (2007). Minimal penalties for Gaussian model selection. *Probab. Theory Related Fields* **138** 33–73. MR2288064

Blake, A. and Zisserman, A. (1987). *Visual Reconstruction*. MIT Press, Cambridge, MA. MR0919733

Boysen, L., Liebscher, V., Munk, A. and Wittich, O. (2007). Scale space consistency of piecewise constant least squares estimators—another look at the regressogram. *IMS Lecture Notes Monograph Ser.* **55** 65–84. IMS, Beachwood, OH.

Braun, J. V., Braun, R. K. and Müller, H.-G. (2000). Multiple change-point fitting via quasilikelihood, with application to DNA sequence segmentation. *Biometrika* **87** 301–314. MR1782480

BURCHARD, H. G. and HALE, D. F. (1975). Piecewise polynomial approximation on optimal meshes. *J. Approximation Theory* **14** 128–147. MR0374761

CHAUDHURI, P. and MARRON, J. S. (2000). Scale space view of curve estimation. *Ann. Statist.* **28** 408–428. MR1790003

CHRISTENSEN, J. and RUDEMO, M. (1996). Multiple change-point analysis of disease incidence rates. *Prev. Vet. Med.* 54–76.

CHU, C., GLAD, I., GODTLIEBSEN, F. and MARRON, J. (1998). Edge-preserving smoothers for image processing. *J. Amer. Statist. Assoc.* **93** 526–541. MR1631321

DAL MASO, G. (1993). *An Introduction to Γ-convergence*. Birkhäuser, Boston. MR1201152

DAVIES, P. L. and KOVAC, A. (2001). Local extremes, runs, strings and multiresolution. *Ann. Statist.* **29** 1–65. MR1833958

DEVORE, R. A. (1998). Nonlinear approximation. In *Acta Numerica 1998. Acta Numer.* **7** 51–150. Cambridge Univ. Press, Cambridge. MR1689432

DEVORE, R. A. and LORENTZ, G. G. (1993). *Constructive Approximation*. Springer, Berlin. MR1261635

DONOHO, D. (2006a). For most large underdetermined systems of equations, the minimal $\ell_1$-norm near-solution approximates the sparsest near-solution. *Comm. Pure Appl. Math.* **59** 907–934. MR2222440

DONOHO, D. (2006b). For most large underdetermined systems of equations, the minimal $\ell_1$-norm solution is the sparsest solution. *Comm. Pure Appl. Math.* **59** 797–829. MR2217606

DONOHO, D. L. (1997). CART and best-ortho-basis: A connection. *Ann. Statist.* **25** 1870–1911. MR1474073

DONOHO, D. L. (1999). Wedgelets: Nearly minimax estimation of edges. *Ann. Statist.* **27** 859–897. MR1724034

DONOHO, D. L. and JOHNSTONE, I. M. (1994). Ideal spatial adaptation by wavelet shrinkage. *Biometrika* **81** 425–455. MR1311089

DONOHO, D. L., JOHNSTONE, I. M., KERKYACHARIAN, G. and PICARD, D. (1995). Wavelet shrinkage: Asymptopia? *J. Roy. Statist. Soc. Ser. B* **57** 301–369. MR1323344

EUBANK, R. L. (1999). *Nonparametric Regression and Spline Smoothing*, 2nd ed. Dekker, New York. MR1680784

FREDKIN, D. and RICE, J. (1992). Baysian restoration and single-channel patch clamp recordings. *Biometrics* **48** 427–428.

FRIEDRICH, F. (2005). Complexity penalized segmentations in 2D. Ph.D. thesis, Institut für Biomathematik und Biometrie an der Gesellschaft für Umwelt und Gesundheit, München-Neuherberg.

FRIEDRICH, F., KEMPE, A., LIEBSCHER, V. and WINKLER, G. (2008). Complexity penalized $m$-estimation: Fast computation. *J. Comput. Graph. Statist.* **17** 1–24.

FÜHR, H., DEMARET, L. and FRIEDRICH, F. (2006). Beyond wavelets: New image representation paradigms. In *Document and Image Compression* (M. Barni, ed.) Chapter 7 179–206. CRC Press.

GEMAN, S. and GEMAN, D. (1984). Stochastic relaxation, Gibbs distributions, and the Bayesian restoration of images. *IEEE Trans. Pattern Anal. Mach. Intell.* **6** 721–741.

GODTLIEBSEN, F., SPJØTVOLL, E. and MARRON, J. S. (1997). A nonlinear Gaussian filter applied to images with discontinuities. *J. Nonparametr. Statist.* **8** 21–43. MR1658111

HALL, P. and TITTERINGTON, D. M. (1992). Edge-preserving and peak-preserving smoothing. *Technometrics* **34** 429–440. MR1190262

HAMPEL, F. R., RONCHETTI, E. M., ROUSSEEUW, P. J. and STAHEL, W. A. (1986). *Robust Statistics*. Wiley, New York. MR0829458

Hess, C. (1996). Epi-convergence of sequences of normal integrands and strong consistency of the maximum likelihood estimator. *Ann. Statist.* **24** 1298–1315. MR1401851

Hinkley, D. V. (1970). Inference about the change-point in a sequence of random variables. *Biometrika* **57** 1–17. MR0273727

Ising, E. (1925). Beitrag zur theorie des ferromagnetismus. *Z. Phys.* **31** 253.

Kohler, M. (1999). Nonparametric estimation of piecewise smooth regression functions. *Statist. Probab. Lett.* **43** 49–55. MR1707251

Künsch, H. R. (1994). Robust priors for smoothing and image restoration. *Ann. Inst. Statist. Math.* **46** 1–19. MR1272743

Loader, C. R. (1996). Change point estimation using nonparametric regression. *Ann. Statist.* **24** 1667–1678. MR1416655

Mammen, E. and van de Geer, S. (1997). Locally adaptive regression splines. *Ann. Statist.* **25** 387–413. MR1429931

Müller, H.-G. (1992). Change-points in nonparametric regression analysis. *Ann. Statist.* **20** 737–761. MR1165590

Müller, H.-G. and Stadtmüller, U. (1999). Discontinuous versus smooth regression. *Ann. Statist.* **27** 299–337. MR1701113

Petrov, V. V. (1975). *Sums of Independent Random Variables*. Springer, New York. MR0388499

Polzehl, J. and Spokoiny, V. (2003). Image denoising: Pointwise adaptive approach. *Ann. Statist.* **31** 30–57. MR1962499

Pötscher, B. and Leeb, H. (2008). Sparse estimators and the oracle property, or the return of Hodges' estimator. *J. Econometrics* **142** 201–211. MR2394290

Potts, R. (1952). Some generalized order-disorder transitions. *Proc. Camb. Philos. Soc.* **48** 106–109. MR0047571

Shao, Q. M. (1995). On a conjecture of Révész. *Proc. Amer. Math. Soc.* **123** 575–582. MR1231304

Spokoiny, V. G. (1998). Estimation of a function with discontinuities via local polynomial fit with an adaptive window choice. *Ann. Statist.* **26** 1356–1378. MR1647669

Tibshirani, R. (1996). Regression shrinkage and selection via the Lasso. *J. Roy. Statist. Soc. Ser. B* **58** 267–288. MR1379242

Tomkins, R. J. (1974). On the law of the iterated logarithm for double sequences of random variables. *Z. Wahrsch. Verw. Gebiete* **30** 303–314. MR0405555

Tukey, J. W. (1961). Curves as parameters, and touch estimation. *Proc. 4th Berkeley Sympos. Math. Statist. and Probab.* **I** 681–694. Univ. California Press, Berkeley. MR0132677

van de Geer, S. (2001). Least squares estimation with complexity penalties. *Math. Methods Statist.* **10** 355–374. MR1867165

Winkler, G. and Liebscher, V. (2002). Smoothers for discontinuous signals. *J. Nonparametr. Statist.* **14** 203–222. MR1905594

Winkler, G., Wittich, O., Liebscher, V. and Kempe, A. (2005). Don't shed tears over breaks. *Jahresber. Deutsch. Math.-Verein.* **107** 57–87. MR2156103

Yao, Y.-C. (1988). Estimating the number of change-points via Schwarz' criterion. *Statist. Probab. Lett.* **6** 181–189. MR0919373

Yao, Y.-C. and Au, S. T. (1989). Least-squares estimation of a step function. *Sankhyā Ser. A* **51** 370–381. MR1175613

L. Boysen
A. Munk
Institute for Mathematical Statistics
Georg–August–Universität Göttingen
Maschmühlenweg 8–10
37073 Göttingen
Germany
E-mail: boysen@math.uni-goettingen.de
munk@math.uni-goettingen.de

A. Kempe
Institute of Biomathematics and Biometry
GSF-National Research Centre
for Environment
Ingolstädter Landstrasse 1
85764 Neuherberg
Germany
E-mail: kempe@gsf.de

V. Liebscher
Departement of Mathematics
and Computer Science
Universität Greifswald
Jahnstrasse 15a
17487 Greifswald
Germany
E-mail: volkmar.liebscher@uni-greifswald.de

O. Wittich
Department of Mathematics
and Computer Science
Technische Universiteit Eindhoven
Den Dolech 2
5600 MB Eindhoven
The Netherlands
E-mail: O.Wittich@tue.nl